\documentclass[11pt]{article}
\usepackage{amssymb}
\usepackage[usenames,dvipsnames]{xcolor}
\usepackage{amsthm}
\usepackage{amsmath}
\usepackage{amsfonts}

\usepackage{graphicx}

 \theoremstyle{definition}

 \numberwithin{equation}{section}

\newtheorem{theorem}{Theorem}[section]

\theoremstyle{definition}
\newtheorem{definition}[theorem]{Definition}

\theoremstyle{remark}
\newtheorem{remark}[theorem]{Remark}

\begin{document}
\title{On steady subsonic flows for Euler-Poisson models}
\author{Shangkun WENG\footnote{Department of mathematics, Harvard University. One, Oxford Street, Cambridge, MA, 02138, USA. {\it Email:
 skwengmath@gmail.com}}}
\date{Harvard University}
\maketitle

\begin{abstract}

In this paper, we are concerned with the structural stability of some steady subsonic solutions for Euler-Poisson system. A steady subsonic solution with subsonic
background charge is proven to be structurally stable with respect to small perturbations of the background charge, the incoming flow angles, the normal electric field and the Bernoulli's function at the inlet and the end pressure at the exit, provided the background solution has a low Mach number and a small electric field. Following the idea developed in \cite{W3}, we give a new formulation for Euler-Poisson equations, which employ the Bernoulli's law to reduce the dimension of the velocity field. The new ingredient in our mathematical analysis is the solvability of a new second
order elliptic system supplemented with oblique derivative conditions at the inlet and Dirichlet boundary conditions at the exit of the nozzle.
\end{abstract}

\begin{center}
\begin{minipage}{5.5in}
2000 Mathematics Subject Classification:35Q35; 35J25;76H05.

\

Key words: Euler-Poisson system, Subsonic flow, second order elliptic system,
Hyperbolic-elliptic coupled system.
\end{minipage}
\end{center}

\

\everymath{\displaystyle}
\newcommand {\eqdef }{\ensuremath {\stackrel {\mathrm {\Delta}}{=}}}

\def\Xint #1{\mathchoice
{\XXint \displaystyle \textstyle {#1}} %
{\XXint \textstyle \scriptstyle {#1}} %
{\XXint \scriptstyle \scriptscriptstyle {#1}} %
{\XXint \scriptscriptstyle \scriptscriptstyle {#1}} %
\!\int}
\def\XXint #1#2#3{{\setbox 0=\hbox {$#1{#2#3}{\int }$}
\vcenter {\hbox {$#2#3$}}\kern -.5\wd 0}}
\def\ddashint {\Xint =}
\def\dashint {\Xint -}
\def\clockint {\Xint \circlearrowright } 
\def\counterint {\Xint \rotcirclearrowleft } 
\def\rotcirclearrowleft {\mathpalette {\RotLSymbol { -30}}\circlearrowleft }
\def\RotLSymbol #1#2#3{\rotatebox [ origin =c ]{#1}{$#2#3$}}

\def\aint{\dashint}

\def\arraystretch{2}
\def\eps{\varepsilon}

\def\s#1{\mathbb{#1}} 
\def\t#1{\tilde{#1}} 
\def\b#1{\overline{#1}}
\def\N{\mathcal{N}} 
\def\M{\mathcal{M}} 
\def\R{{\mathbb{R}}}

\def\ba{\begin{array}}
\def\ea{\end{array}}
\def\be{\begin{equation}}
\def\ee{\end{equation}}

\def\bes{\begin{mysubequations}}
\def\ees{\end{mysubequations}}

\def\cz#1{\|#1\|_{C^{0,\alpha}}}
\def\ca#1{\|#1\|_{C^{1,\alpha}}}
\def\cb#1{\|#1\|_{C^{2,\alpha}}}

\def\lb#1{\|#1\|_{L^2}}
\def\ha#1{\|#1\|_{H^1}}
\def\hb#1{\|#1\|_{H^2}}

\def\cin{\subset\subset}
\def\Ld{\Lambda}
\def\ld{\lambda}
\def\ol{{\Omega_L}}
\def\sla{{S_L^-}}
\def\slb{{S_L^+}}
\def\e{\varepsilon}
\def\C{\mathbf{C}} 
\def\cl#1{\overline{#1}}
\def\ra{\rightarrow}
\def\xra{\xrightarrow}
\def\g{\nabla}
\def\a{\alpha}
\def\b{\beta}
\def\d{\delta}
\def\th{\theta}
\def\fai{\varphi}
\def\O{\Omega}
\def\f{\frac}
\def\p{\partial}
\def\disp{\displaystyle}

\def\H{\Theta} 


\section{Introduction}\label{EulerPoissonintroduction}\hspace*{\parindent}

The Euler-Poisson (or hydrodynamical) model for a unipolar semiconductor in the isentropic steady state case reads
 \be\label{ea1} \left\{\ba{l}
 (\rho u_1)_{x_1}+(\rho u_2)_{x_2}+(\rho u_3)_{x_3}=0,\\
 (\rho u_1^2)_{x_1}+(\rho u_1u_2)_{x_2}+(\rho u_1u_3)_{x_3}+p_{x_1}=\rho \varphi_{x_1},\\
 (\rho u_1u_2)_{x_1}+(\rho u_2^2)_{x_2}+(\rho u_2u_3)_{x_3}+p_{x_2}=\rho \varphi_{x_2},\\
 (\rho u_1u_3)_{x_1}+(\rho u_2u_3)_{x_2}+(\rho u_3^2)_{x_3}+p_{x_3}=\rho \varphi_{x_3},\\
 \Delta \varphi= \rho-b(x).
\ea\right.
\ee where $\rho, {\bf u},\varphi$ denotes the electron density, the average electron velocity and the electrostatic
potential, $b(x)$ be the prescribed ion background density (doping profile). The electric field ${\bf E}$ is given
by ${\bf E}=\nabla \varphi$. The energy equation of the hydrodynamic
model is replaced by the pressure-density relation $p=p(\rho)$. One may refer to \cite{M1,MRS} for more details on physical background.

In this paper, we are concerned with the structural stability of some steady subsonic solutions for Euler-Poisson system. Indeed, we will construct a subsonic solution to (\ref{ea1}) in a rectangular cylinder by imposing suitable
boundary conditions, which is also required to be close to
some special subsonic solutions. The boundary conditions we will impose have their origins on our previous works in \cite{DWX,W1,W2,W3} and another important work \cite{LX} on transonic shock solutions for (\ref{ea1}). As we have done in \cite{W3},
we will prescribe the incoming flow angles $\beta_2=\f{u_2}{u_1},\beta_3=\f{u_3}{u_1}$ and Bernoulli's function $B=\f{1}{2}(u_1^2+u_2^2+u_3^2)+h(\rho)-\varphi$ at the inlet, the natural slip boundary conditions on the nozzle walls and the end pressure at the exit.
For the electric field ${\bf E}=\nabla \varphi$, we prescribe the normal component of ${\bf E}$ at the inlet and nozzle walls,
and $\varphi$ at the exit of nozzle.

The one-dimensional steady state isentropic flow hydrodynamical model was analyzed in \cite{DM1} and the three-dimensional irrotational case was disucssed in \cite{DM2}. In both papers, existence and uniqueness results for small data generating subsonic flow were proven. In \cite{M2}, the author obtained the 2-D existence of smooth solution under some smallness assumptions on the prescribed outflow current and the gradient of the velocity relaxation time. Also, the author used a smallness assumption on the physical parameter multiplying the drift-term in the velocity equations instead in the irrotational 3-D case and obtained similar results. For the full hydrodynamic model,Yeh \cite{Yeh2} showed the existence of a unique strong solution in
several space dimensions if the flow is subsonic, the ambient temperature is large enough, and the vorticity on the inflow boundary and the
variation of the electron density on the boundary are sufficiently small. Zhu and Hattori \cite{ZH} proved the existence of classical subsonic solutions in one space dimension for the whole space problem under the additional assumption that the doping profile be close to a constant. See also \cite{ABJM} and the reference therein for more details.

Recently, there are some important progress on transonic shock solutions for 1-D Euler-Poisson systems \cite{LX,LRXX}. Gamba \cite{Gamba} showed existence of steady-state solutions in the transonic case by means of the vanishing viscosity method. However, the solutions as the limit of vanishing
viscosity may contain boundary layers and more than one transonic shock.
In \cite{LX}, Luo and Xin gave a thorough study of the structure of the solutions to boundary value problems for 1-D Euler-Poisson system
for different situations when the density of fixed, positively
charged background ions is in supersonic and subsonic regimes. The existence,
non-existence, uniqueness and non-uniqueness of solutions with transonic shocks were obtained
according to the different cases of boundary data and physical interval length. The solutions they constructed contained exactly one transonic shock in the interval $[0,L]$. Moreover,
they can determined the shock location by the boundary data and $L$. In \cite{LRXX}, the authors investigated structural and dynamical stabilities of steady transonic shock
solutions for one-dimensional Euler-Poisson systems. It was shown that a
steady transonic shock solution with a supersonic background charge
was structurally stable with respect to small perturbations of the background
charge, provided that the electric field is positive at the shock location. Furthermore, any
steady transonic shock solution with a supersonic background charge was proved to
be dynamically and exponentially stable with respect to small perturbations of the
initial data, provided the electric field is not too negative at the shock location.
The remaining natural question is the structural stability of these transonic shock solutions in a multi-dimensional domain supplemented with
suitable boundary conditions, a considerably more difficult matter. We start to investigate the structural stability of some special subsonic solution to take a close look at the physically acceptable boundary conditions we
should prescribe and hope this may serve as a building block toward that challenging goal.

In \cite{DWX,W1,W2}, we have characterized a set of physically acceptable boundary conditions that ensure the existence and uniqueness of a subsonic irrotational flow in a finitely long flat nozzle by prescribing the incoming flow angle and the Bernoulli's function at the inlet and the end pressure at the exit. In \cite{W3}, we have developed a new formulation for the three dimensional Euler equations. The key idea in our formulation is to use the Bernoulli's law to reduce the dimension of the velocity field by defining new variables $(1,\beta_2=\frac{u_2}{u_1},\beta_3=\frac{u_3}{u_1})$ and replacing $u_1$ by the Bernoulli's function $B$. We find a conserved quantity for flows with a constant Bernoulli's function, which behaves like the scaled vorticity in the 2-D case. Moreover, a system of new conservation laws can be derived, which is new even in the two dimensional case. Following the ideas developed in \cite{W3}, we reformulate the Euler-Poisson equations in terms of $(s=\ln \rho, \beta_2,\beta_3, B,\varphi)$ and replace $u_1$ by $B$ through $u_1^2=\frac{2(B+\varphi-h(\rho))}{1+\beta_2^2+\beta_3^2}$. Here one should note that due to the Bernoulli's law (\ref{eb00}), the Bernoulli's function $B$ will possess the same regularity as the boundary data at the inlet, which is quite different from the velocity field.  In this way, we can explore the role of the Bernoulli's law in greater depth and hope that may simplify the Euler equations a little bit. We can also find a new quantity $W=\p_2\beta_3-\p_3\beta_2+\beta_3\p_1\beta_2-\beta_2\p_1\beta_3$, which is conserved along the particle path at least for flows with a constant Bernoulli's function (see (\ref{eb26})). In subsonic region, we find that the density and the potential function are coupled together to satisfy an second order elliptic system, which will possess good regularity. Roughly speaking, the equation (\ref{eb26}) shows that $W$ possess one order lower regularity than those of the density and the potential function. This, together with the first three equations in (\ref{eb22}) form an divergent-curl system for $(\beta_2,\beta_3)$, which indicates the velocity field $(\beta_2,\beta_3)$ should possess the same regularity as those of the density and the potential function.

We make some comments on our proof. The new ingredient in our mathematical analysis is the solvability of some second
order elliptic system supplemented with oblique derivative conditions at the inlet and Dirichlet boundary conditions at the exit of the nozzle.
Indeed, both $\rho$ and $\fai$ satisfy some second order elliptic equations and they are coupled together to form a second order elliptic system.
Besides this, the boundary conditions for $\rho$ and $\fai$ are also coupled together at the inlet of the nozzle. This elliptic system is neither weakly-coupled nor cooperative in the sense of \cite{deFM}. To guarantee the uniqueness, some smallness assumptions are prescribed on the background solution. That means our background solution should have a low Mach number and a small electric field.

This paper proceeds as follows. The second section contains a new formulation for Euler-Poisson system and a statement of our background solution.
The last section focus on the structural stability of our background solution with respect to small perturbations of the background charge, the incoming flow angles and the end pressure,
provided the background solution has a low Mach number and a small electric field.

\section{Preliminary}\label{EulerPoissonprelimiary}\hspace*{\parindent}

\subsection{A new formulation for Euler-Poisson equations}\label{NewFormulation} \hspace*{\parindent}

Following \cite{W3}, we develop a new formulation for Euler-Poisson equations, which employ the Bernoulli's law to reduce the dimension of the velocity field. By employing the following identity in vector calculus
$${\bf u}\cdot \nabla {\bf u}=\nabla (\frac{1}{2}|{\bf u}|^2)-{\bf u}\times curl {\bf u},$$
and the momentum equation in (\ref{ea1}), we have
\be\label{eb0}
\nabla (\frac{1}{2}|{\bf u}|^2+h(\rho)-\varphi)-{\bf u}\times curl {\bf u}=0.
\ee
Then the Bernoulli's law holds:
\begin{equation}\label{eb00}
{\bf u}\cdot \nabla B=0.
\end{equation}
Here the Bernoulli's function is defined to be $B=\f{1}{2}(u_1^2+u_2^2+u_3^2)-\fai+h(\rho)$,
where $h'(\rho)=\f{c^2(\rho)}{\rho}=\f{p'(\rho)}{\rho}$.

Define the following three new variables: $\beta_2=\f{u_2}{u_1},\beta_3=\f{u_3}{u_1}, s=\ln \rho$. Then $u_1$ has a simple expression: $u_1^2=\f{2(B+\varphi-h(\rho))}{1+\beta_2^2+\beta_3^2}$.

Multiplying the first equation  in (\ref{ea1}) by $\f{-u_{i-1}}{\rho
u_1^2}$, dividing the i-th equation in (\ref{ea1}) by $\rho u_1^2$ and adding them together for $i=2,3,4$, we obtain the following new system:
\be\label{eb22} \left\{\ba{l}
\p_1s+\beta_2\p_2s+\beta_3\p_3s-\f{c^2(\rho)}{u_1^2}\p_1s+\p_2\beta_2+\p_3\beta_3=-\f{1}{u_1^2}\p_1\varphi,\\
\p_1\beta_2+\beta_2\p_2\beta_2+\beta_3\p_3\beta_2-\f{c^2(\rho)}{u_1^2}\beta_2\p_1s+\f{c^2(\rho)}{u_1^2}\p_2s+\f{\beta_2}{u_1^2}\p_1\varphi-\f{1}{u_1^2}\p_2\varphi=0,\\
\p_1\beta_3+\beta_2\p_2\beta_3+\beta_3\p_3\beta_3-\f{c^2(\rho)}{u_1^2}\beta_3\p_1s+\f{c^2(\rho)}{u_1^2}\p_3s+\f{\beta_3}{u_1^2}\p_1\varphi-\f{1}{u_1^2}\p_3\varphi=0,\\
\p_1B+\beta_2\p_2B+\beta_3\p_3B=0. \\
\Delta \varphi=\rho-b(x).
\ea\right.
\ee

Then it is easy to show $s$ and $\varphi$ will satisfy the following second order elliptic system in the subsonic region, i.e. $|{\bf u}|\leq c(\rho)$,
\be\label{eb24} \left\{\ba{l}
\p_1\bigg((\f{c^2}{u_1^2}-1)\p_1 s-\beta_2\p_2 s-\beta_3\p_3 s\bigg)+\p_2\bigg(-\beta_2\p_1 s+(\f{c^2(\rho)}{u_1^2}-\beta_2^2)\p_2 s-\beta_2\beta_3\p_3 s\bigg)\\
+\p_3\bigg(-\beta_3\p_1 s-\beta_2\beta_3\p_2 s+(\f{c^2(\rho)}{u_1^2}-\beta_3^2)\p_3 s\bigg)-\bigg((\f{c^2}{u_1^2}-1)\p_1 s-\beta_2\p_2 s-\beta_3\p_3 s\bigg)^2
\\+\bigg((\p_2\beta_2)^2+(\p_3\beta_3)^2+2\p_2\beta_3\p_3\beta_2\bigg)-\f{1}{u_1^2}({e}^s-b(x))-\nabla \varphi\times\nabla \f{1}{u_1^2}=0,\\
\Delta \varphi={e}^s-b(x).
\ea\right.\ee

An important observation made in \cite{W3} is the quantity $W=\p_2\beta_3-\p_3\beta_2+\beta_3\p_1\beta_2-\beta_2\p_1\beta_3$ is conserved along the particle path for flow with a constant Bernoulli's function. To simplify
the notation, we set $G=1+\beta_2^2+\beta_3^2$ and ${\bf D}= \p_1+\beta_2\p_2+\beta_3\p_3$. Plugging $u_1^2=\f{2(B+\varphi-h(\rho))}{1+\beta_2^2+\beta_3^2}$ into the second and third equation in (\ref{eb22}), then we obtain
\be\label{eb25} \left\{\ba{l}
{\bf D}\beta_2-\f{c^2(\rho)}{2(B+\varphi-h(\rho))}G(\beta_2\p_1s-\p_2s)+\f{G}{2(B+\varphi-h(\rho))}(\beta_2\p_1\varphi-\p_2\varphi)=0,\\
{\bf D}\beta_3-\f{c^2(\rho)}{2(B+\varphi-h(\rho))}G(\beta_3\p_1s-\p_3s)+\f{G}{2(B+\varphi-h(\rho))}(\beta_3\p_1\varphi-\p_3\varphi)=0.
\ea\right.\ee

Apply $\beta_3\p_1-\p_3$ and $-\beta_2\p_1+\p_2$ to the above two
equations respectively and add them together, one can
show that $W$ satisfies the following equation.
\be\label{eb26}{\bf D}\bigg(\f{W}{\rho G}\bigg)+\f{1}{2\rho(B+\varphi-h(\rho))^2}\bigg((1,\beta_2,\beta_3)\cdot[(\f{c^2(\rho)}{\rho}\nabla \rho-\nabla\varphi)\times\nabla B]\bigg)=0.\ee

The verification of (\ref{eb26}) is similar to the calculation in the appendix in \cite{W3}. For completeness, we also give a detailed calculation for (\ref{eb26}) in the appendix.

Suppose that $B\equiv const$, then (\ref{eb26}) reduces to
\be\label{eb261}
{\bf D}\bigg(\f{W}{\rho G}\bigg)=0.
\ee
This implies $\f{W}{\rho G}$ is conserved along the particle path. Indeed, since $B\equiv const$, by (\ref{eb0}), the vorticity field is parallel to the velocity field. So we may assume that there exists a real function $\mu(x)$ such that $curl{\bf u}= \mu(x){\bf u}$. By simple calculations, we have
\begin{eqnarray}\label{eb262}
\f{W}{\rho G}&=& \f{curl{\bf u}\cdot{\bf u}}{\rho u_1^2 G}
=\f{\mu(x)|{\bf u}|^2}{\rho u_1^2 G}=\f{\mu(x)}{\rho}.
\end{eqnarray}

Suppose that $B\equiv const$, then (\ref{eb26}) reduces to
\be\label{eb263}
{\bf D}\bigg(\f{\mu(x)}{\rho}\bigg)=0.
\ee
Historically, the stationary solution of Euler equations with the vorticity field paralleling to the velocity field was called Beltrami flow and have been investigated for over a century. One may refer to Arnol'd \cite{Arnol'd}, Constantin and Majda \cite{CM} for more details. Here we emphasize that our calculations works for any general Euler-Poisson flows.

\begin{remark}{\it As we have discussed in the introduction, $W$ may help to raise the regularity of the velocity field $\beta_2$ and $\beta_3$. In subsonic region, the equation (\ref{eb24}) shows that the density and the potential function are coupled together to satisfy an second order elliptic system, which will possess good regularity. Roughly speaking, the equation (\ref{eb26}) shows that $W$ possess one order lower regularity than those of the density and the potential function. This, together with the first three equations in (\ref{eb22}) form an divergent-curl system for $(\beta_2,\beta_3)$, which indicates the velocity field $(\beta_2,\beta_3)$ should possess the same regularity as those of the density and the potential function.}
\end{remark}

\begin{remark}{\it
It is easy to see that $G$ satisfies the following
Riccati-type equation:
\be\label{eb27}
{\bf D}G-\f{1}{(B+\varphi-h(\rho))}G^2\p_1(h(\rho)-\varphi)-\f{1}{(B+\varphi-h(\rho))}G {\bf D}(h(\rho)-\varphi)=0.
\ee
One may expect some blow-up results. However, the ambiguous sign of $\p_1 (h(\rho)-\varphi)$ in the coefficient of $G^2$ in (\ref{eb27}) makes the whole argument nontrivial. One should note that $G$ blows up means that the fluid will turn around in the flow region.}
\end{remark}

\begin{remark}{\it
In \cite{W3}, we can also derive a system of new conservation laws for compressible Euler equations. However, we can not obtain similar conservation laws for the Euler-Poisson equations due to the effect of the electrostatic potential.}
\end{remark}

\begin{remark}
{\it $\rho$ and $\fai$ are coupled together to form a second order elliptic system. This elliptic system is neither weakly-coupled
nor cooperative in the sense of \cite{deFM}. Hence the uniqueness for (\ref{eb24}) with suitable boundary conditions will be a main obstacle in our analysis.
}
\end{remark}

\subsection{Background solutions}\label{EulerPoissonbackgroundsolutions}\hspace*{\parindent}

In \cite{LX}, the author consider the initial value problem for the following 1-D Euler-Poisson equations
\begin{equation}\label{SteadyEP}
\begin{cases}
& (\rho u)_{x}=0, \\
& (p(\rho )+\rho u^{2})_{x}=\rho E, \\
& E_{x}=\rho-b_0.\\
& (\rho ,u,E)(0)=(\rho _{I},u_{I},E_I).%
\end{cases}
\end{equation}%
Here $b_0$ is a positive constant.

Assume $u_{I}>0$. By the first equation in (\ref{SteadyEP}), we have $\rho u(x)=J=\rho _{I}u_{I}$ and the velocity is given by
\begin{equation}\label{velocity}
u=J/\rho .
\end{equation}
Thus the boundary value problem  (\ref{SteadyEP}) reduces to
\begin{equation}\label{SimpleSteadyEP}
\begin{cases}
& (p(\rho )+\frac{J^{2}}{\rho })_{x}=\rho E, \\
& E_{x}=\rho -b_0,\\
&(\rho ,E)(0)=(\rho _{I}, E_I ).%
\end{cases}
\end{equation}%

Use the terminology from gas dynamics to call $c=\sqrt{p^{\prime }(\rho )}
$ the sound speed. There is a unique solution $\rho =\rho_{s}$ for the
equation
\begin{equation} \label{sonic}
p^{\prime }(\rho )=J^{2}/\rho^2,
\end{equation}%
which is the sonic state (recall that $J=\rho u$).  Later on, the flow is called supersonic if
\begin{equation}
p^{\prime }(\rho )<J^{2}/\rho^2,\ i.e.\ \rho <\rho_{s}.
\end{equation}%
Similarly, if
\begin{equation}
p^{\prime }(\rho )>J^{2}/\rho^2,\ i.e. \ \rho >\rho_{s},
\end{equation}%
then the flow is said to be subsonic.

The solution of (\ref{SimpleSteadyEP}) was analyzed in $(\rho,E)-$phase plane. Any trajectory in $(\rho,E)-$plane satisfies the following equation,
\begin{equation}
d(\f{1}{2}E^2-H(\rho))=0, \text{where}\ \ H'(\rho)=\f{\rho-b_0}{\rho}(p'(\rho-\f{J^2}{\rho^2})).
\end{equation}

The trajectory passing through the point $(\rho_I,E_I)$ with $\rho_I>0$ is given by
\be
\f{1}{2}E^2-\int_{\rho_I}^{\rho}H'(s)ds=\f{1}{2}E_I^2.
\ee
We only consider the case $b_0>\rho_s$, i.e. $b_0$ is in subsonic region.

\begin{definition}
The critical trajectory (for the case $b_0>\rho_s$) is the trajectory passing through the point $(b_0,0)$ with the equation:
\be
\f{1}{2}E^2-\int_{b_0}^{\rho}H'(s)ds=0.
\ee
\end{definition}
We are concerned with the following two cases, that is case c1) and c4) in the subsection 2.2 in \cite{LX}.
Suppose $(\rho_I,E_I)$ is on the critical supersonic trajectory, i.e.
$$\f{1}{2}E_I^2+\int_{\rho_I}^{b_0}H'(s)ds=0.$$
\begin{itemize}
\item  $\rho_s<\rho_I<b_0, E_I>0.$\\

In this case, initial value problem (\ref{SimpleSteadyEP}) admits a unique subsonic solution $(\rho,E)$ for all $x\geq 0$. Moreover,
\be
\rho_x>0, E_x<0, x>0,\\
\lim_{x\rightarrow \infty}(\rho,E)(x)=(b_0,0).
\ee

\item $\rho_I>b_0,E_I<0.$\\

In this case, initial value problem (\ref{SimpleSteadyEP}) admits a unique subsonic solution $(\rho,E)$ for all $x\geq 0$. Moreover,
\be
\rho_x>0, E_x>0, x>0,\\
\lim_{x\rightarrow \infty}(\rho,E)(x)=(b_0,0).
\ee
\end{itemize}

Due to our technical reasons, we need to choose some special background solutions. Given $b_0>0$, we can choose $J>0$ such that $b_0>\rho_s$, i.e. $b_0$ is in subsonic region.
We take $\rho_I$ be close to $b_0$, so that $0<u_I<\f{1}{3}\sqrt{\f{1}{b_0C_{\O_e}}}$(we can do it by choose a small $J$). Then we choose
$0<|E_I|<\f{1}{C(b_0)C_{\O_e}^3}$ and $(\rho_I,E_I)$ is on the critical supersonic trajectory. By the above analysis, we can find a unique subsonic solution $(\rho_0,E_0)$
to (\ref{SimpleSteadyEP}) for all $x\geq 0$. Moreover, we have $\lim_{x\rightarrow \infty}(\rho,E)(x)=(b_0,0)$. Restricted to $[0,1]$, $(\rho_0,E_0)$ will satisfy the following
properties:
\be\label{eb25}
\begin{cases}
\f{b_0}{2}<\min_{x\in[0,1]}\rho_0(x)<\max_{x\in[0,1]}\rho_0(x)<\f{3}{2}b_0.\\
\f{1}{6}\sqrt{\f{1}{b_0C_{\O_e}}}<\min_{x\in[0,1]}u_0(x)<\max_{x\in[0,1]}u_0(x)<\f{2}{3}\sqrt{\f{1}{b_0C_{\O_e}}}.\\
\max_{x\in[0,1]}E_0(x)<\f{1}{C(b_0)C_{\O_e}^3}.
\end{cases}
\ee
Here $C(b_0)$ is a smooth function of $b_0$. The constant $C_{\O_e}$ is the least number such that the following inequality holds:
$$
\|U(x)\|_{L^2(\O_e)}^2\leq C_{\O_e}\|\nabla U(x)\|_{L^2(\O_e)}^2,\forall \ U(x)\in\mathbb{H}=\{U(x)\in H^1(\O_e):U(1,x_2,x_3)=0.\}.
$$

For convenience, we introduce the electrostatic potential $\varphi_0(x)$ satisfying $\varphi_0'(x)=E_0(x)$ with $\varphi(1)=0$.
\begin{remark}
{\it
The background solution we have chosen has a low Mach number and a small electric field in the sense of (\ref{eb25}).
}
\end{remark}
\section{Structural stability of background solutions}\label{EulerPoissonstructuralstability} \hspace*{\parindent}

In this section, we are concerned with the structural stability of our background solutions for Euler-Poisson system.
Indeed, we will construct a subsonic solution to (\ref{ea1}) in a rectangular cylinder by imposing suitable
boundary conditions at the inlet and exit, which is also required to be close to
our background solutions $(\rho_0(x_1),u_0(x_1),\fai_0(x_1))$. The rectangular cylinder will be $\O=[0,1]\times[0,1]\times[0,1]$.
We also set $B_0=\f{1}{2}u_0^2(x_1)+h(\rho_0(x_1))-\fai_0(x_1)$ is a constant, $s_0(x_1)=\ln\rho_0(x_1)$.

At the inlet of the nozzle $x_1=0$, we impose the flow angles and the Bernoulli's function:
\be\label{ee1}\left\{\ba{l}
\beta_i(0,x_2,x_3)=\epsilon \beta_{i}^{in}(x_2,x_3), i=2,3,\\
B(0,x_2,x_3)= B_0+\epsilon B^{in}(x_2,x_3).
\ea\right.\ee
Here the compatibility conditions should be satisfied:
\be\label{ee100}\left\{\ba{l}
\p_2^j\beta_{2}^{in}(0,x_3)=\p_2^j\beta_{2}^{in}(1,x_3)=0, \\
\p_3^j\beta_{3}^{in}(x_2,0)=\p_3^j\beta_{3}^{in}(x_2,1)=0,\ \ j=0,2,\\
\p_2^k B^{in}(0,x_3)=\p_2^k B^{in}(1,x_3)=\p_3^k B^{in}(x_2,0)=\p_3^k B^{in}(x_2,1)=0,\ \ k=1,3.
\ea\right.\ee

At the exit of the nozzle $x_1=1$, we prescribe the end pressure:
\be\label{ee2}
p(1,x_2,x_3)=\f{1}{\gamma} e^{\gamma(s_0(1)+\epsilon s_e(x_2,x_3))}.
\ee
Here we also require that $p_e$ satisfies the following compatibility conditions:
\be\label{ee101}\left\{\ba{l}
\p_2^j s_e(0,x_3)=\p_2^j s_e(1,x_3)=0,\\
\p_3^j s_e(x_2,0)=\p_3^j s_e(x_2,1)=0, \ \ j=1,3.
\ea\right.\ee

While on the nozzle walls, the usual slip boundary condition is imposed:
\be\label{ee3}\left\{\ba{l}
u_2(x_1,0,x_3)=u_2(x_1,1,x_3)=0,\\
u_3(x_1,x_2,0)=u_3(x_1,x_2,1)=0.
\ea\right.\ee

For the electric field ${\bf E}=\nabla \fai$, we impose the normal component of ${\bf E}$ at the inlet and the nozzle walls,
 while at the exit, we prescribe the Dirichlet boundary condition for $\fai$:
\be\label{ee31}\left\{\ba{l}
\p_1\fai(0,x_2,x_3)=E_I+\epsilon E^{in}(x_2,x_3),\\
\p_2\fai(x_1,0,x_3)=\p_2\fai(x_1,1,x_3)=0,\\
\p_3\fai(x_1,x_2,0)=\p_3\fai(x_1,x_2,1)=0,\\
\fai(1,x_2,x_3)=0.
\ea\right.\ee
Here $E^{in}(x_2,x_3)$ should satisfy the following compatibility conditions as $B^{in}$:
\be\label{ee32}
\p_2^k E^{in}(0,x_3)=\p_2^k E^{in}(1,x_3)=\p_3^k E^{in}(x_2,0)=\p_3^k E^{in}(x_2,1)=0,\ \ k=1,3.
\ee

The prescribed ion background density $b(x)$ is taken to be $b(x)=b_0+\epsilon \t b(x)$. Here $\t b(x)\in C^{1,\alpha}(\O)$ should satisfy the
following compatibility conditions:
\be\label{ee33}
\p_2 \t b(x_1,0,x_3)=\p_2 \t b(x_1,1,x_3)=\p_3 \t b(x_1,x_2,0)=\p_3 \t b(x_1,x_2,1)=0.
\ee

Mathematically, we are going to prove that (\ref{ea1}) with boundary conditions (\ref{ee1}),(\ref{ee2}), (\ref{ee3}) and (\ref{ee31})
satisfying compatibility conditions (\ref{ee100}), (\ref{ee101}) and (\ref{ee32}), has a unique subsonic solution.

\subsection{Extension to the domain $\O_e=[0,1]\times \mathrm{T}^2$}\label{EulerPoissonextension}\hspace*{\parindent}

Suppose the flow $(\rho, u_1,u_2,u_3)\in C^{3,\alpha}(\bar \O)\times C^{2,\alpha}(\bar \O)^3$ we will construct has the following properties:
\be\label{ee41} \left\{\ba{l}
\p_2^j(\rho,u_1,\fai)(x_1,0,x_3)=\p_2^j(\rho,u_1,\fai)(x_1,1,x_3)=0,\\
\p_3^j(\rho,u_1,\fai)(x_1,x_2,0)=\p_3^j(\rho,u_1,\fai)(x_1,x_2,1)=0,\ \ j=1,3,\\
\p_2^ku_2(x_1,0,x_3)=\p_2^ku_2(x_1,1,x_3)=0,\\
\p_3^ku_3(x_1,x_2,0)=\p_3^ku_3(x_1,x_2,1)=0,\ \ k=0,2.
\ea\right. \ee
Then we may extend $(\rho, u_1,u_2,u_3)$ in the following way, to $(\hat\rho, \hat u_1,\hat u_2,\hat u_3)\in C^{3,\alpha}([0,1]\times \mathbb{R}^2)\times C^{2,\alpha}([0,1]\times \mathbb{R}^2)^3$:

For $(x_1,x_2,x_3)\in [0,1]\times [-1,1]\times[-1,1]$, we define $(\hat\rho,\hat u_1,\hat u_2,\hat u_3)$ as follows
$$
(\hat\rho,\hat u_1,\hat u_2,\hat u_3,\hat \fai)({\bf x})=\left\{
    \ba{ll}
    (\rho,u_1, u_2, u_3,\fai)(x_1,x_2,x_3),\ \ \ \text{ if } (x_2,x_3)\in [0,1]\times [0,1],\\
    (\rho,u_1, -u_2, u_3,\fai)(x_1,-x_2,x_3),\ \ \text{ if } (x_2,x_3)\in [-1,0]\times [0,1],\\
    (\rho,u_1, u_2, -u_3,\fai)(x_1,x_2,-x_3),\ \ \text{ if } (x_2,x_3)\in [0,1]\times [-1,0],\\
    (\rho,u_1, -u_2, -u_3,\fai)(x_1,-x_2,-x_3),\ \text{ if } (x_2,x_3)\in [-1,0]\times [-1,0].
    \ea
    \right.
$$
Then we extend $(\hat\rho,,\hat\fai), \hat u_1,\hat u_2,\hat u_3$ periodically to $[0,1]\times \mathbb{R}^2$ with period $2$.
It is easy to verify that $(\hat\rho,,\hat\fai), \hat u_1,\hat u_2,\hat u_3)$
will belong to $C^{3,\alpha}([0,1]\times \mathbb{R}^2)^2\times C^{2,\alpha}([0,1]\times \mathbb{R}^2)^3$.
Moreover, $(\hat\rho, \hat\fai, \hat u_1,\hat u_2,\hat u_3)$ will also satisfy Euler-Poisson equations.
Due to these reasons, one may directly work on the domain $[0,1]\times \mathrm{T}^2$ (Here $\mathrm{T}^2$ is a 2-torus),
so that the difficulty caused by corner singularity and slip boundary conditions will be avoided.
We may extend $(\beta_{2}^{in},\beta_{3}^{in},B^{in},E^{in}, s_e)$ to $\mathrm{T}^2$,
which will still be denoted by $(\beta_{2}^{in},\beta_{3}^{in},B^{in},E^{in}, s_e)$.

\subsection{Main results}\label{EulerPoissonmainresults}\hspace*{\parindent}

The main result is the following existence and uniqueness theorem.
\begin{theorem}\label{EPTH}
Given $(\beta_{2}^{in},\beta_{3}^{in},B^{in},E^{in}, s_e)\in C^{3,\alpha}(\mathrm{T}^2)$ and $\t b(x)\in C^{1,\alpha}(\O_e)$, there exists a positive small number $\epsilon_0$, which depends on the background subsonic state $(\rho_0,u_0,0,0,\fai_0)$
 and $(\beta_{2}^{in},\beta_{3}^{in},B^{in},E^{in}, s_e,\t b)$, such that if $0<\epsilon<\epsilon_0$,
 then there exists a unique smooth subsonic flow $(u_1,u_2,u_3,\rho,\fai)\in C^{2,\alpha}(\O_e)^3\times C^{3,\alpha}(\O_e)^2$ to (\ref{ea1})
 satisfying boundary conditions (\ref{ee1})- (\ref{ee3}) and (\ref{ee31}). Moreover, the following estimate holds:
\be\label{ee50}
\|(u_1,u_2,u_3)-(u_0,0,0)\|_{C^{2,\alpha}(\O_e)}+\|\rho-\rho_0\|_{C^{3,\alpha}(\O_e)}+\|\fai-\fai_0\|_{C^{3,\alpha}(\O_e)}\leq C\epsilon.
\ee
Here $C$ is a constant depending on $(\rho_0,u_0,0,\fai_0)$ and $(\beta_{2}^{in},\beta_{3}^{in},B^{in},E^{in} ,s_e,\t b)$.
\end{theorem}

Before we start to prove our main theorem, we need to make some preparations.

Define $W_1=s-s_0, W_2=\beta_2, W_3=\beta_3, W_4=B-B_0, W_5=\fai-\fai_0$, we can derive the equations satisfied by $W_i,i=1,\cdots,5$.
\be\label{ee23} \left\{\ba{l}
(1-\f{c^2(\rho_0)}{u_0^2})\p_1W_1+\p_2W_2+\p_3W_3+\f{1}{u_0^2}\p_1W_5+d_1W_1+d_4W_4+d_5W_5=-F_1({\bf  W},\nabla {\bf  W}),\\
\p_1W_2+W_2\p_2W_2+W_3\p_3W_2+\f{c^2(\rho_0)}{u_0^2}\p_2W_1-\f{1}{u_0^2}\p_2W_5+d_2W_2=F_2({\bf  W},\nabla {\bf  W}),\\
\p_1W_3+W_2\p_2W_3+W_3\p_3W_3+\f{c^2(\rho_0)}{u_0^2}\p_3W_1-\f{1}{u_0^2}\p_3W_5+d_3W_3=F_3({\bf  W},\nabla {\bf  W}),\\
\p_1W_4+W_2\p_2W_4+W_3\p_3W_4=0. \\
\Delta W_5=\rho_0 W_1-\epsilon\t b(x)+F_5({\bf  W},\nabla {\bf  W}).
\ea\right. \ee

Here
\be\label{ee24} \left\{\ba{l}
F_1({\bf  W},\nabla {\bf  W})=W_2\p_2W_1+W_3\p_3W_1-(\f{c^2(\rho)}{u_1^2}-\f{c^2(\rho_0)}{u_0^2})\p_1W_1+(\f{1}{u_1^2}-\f{1}{u_0^2})\p_1 W_5\\-(\f{c^2(\rho)}{u_1^2}-\f{c^2(\rho_0)}{u_0^2}-a_2W_1+a_1W_4+a_1W_5)\p_1s_0+(\f{1}{u_1^2}-\f{1}{u_0^2}-b_2W_1+b_1W_4+b_1W_5)\p_1\fai_0,\\
F_2({\bf  W},\nabla {\bf  W})=\f{c^2(\rho)}{u_1^2}W_2\p_1W_1-\f{W_2}{u_1^2}\p_1W_5-(\f{c^2(\rho)}{u_1^2}-\f{c^2(\rho_0)}{u_0^2})\p_2W_1\\+(\f{1}{u_1^2}-\f{1}{u_0^2})\p_2 W_5+(\f{c^2(\rho)}{u_1^2}-\f{c^2(\rho_0)}{u_0^2})\p_1s_0W_2-(\f{1}{u_1^2}-\f{c^2(\rho_0)}{u_0^2})\p_1\fai_0W_2,\\
F_3({\bf  W},\nabla {\bf  W})=\f{c^2(\rho)}{u_1^2}W_3\p_1W_1-\f{W_3}{u_1^2}\p_1W_5-(\f{c^2(\rho)}{u_1^2}-\f{c^2(\rho_0)}{u_0^2})\p_3W_1\\+(\f{1}{u_1^2}-\f{1}{u_0^2})\p_3 W_5+(\f{c^2(\rho)}{u_1^2}-\f{c^2(\rho_0)}{u_0^2})\p_1s_0W_3-(\f{1}{u_1^2}-\f{c^2(\rho_0)}{u_0^2})\p_1\fai_0W_3,\\
F_5({\bf  W},\nabla {\bf  W})=\rho-\rho_0-\rho_0 W_1.
\ea\right.\ee
and
\be\label{ef3} \left\{\ba{l}
a_1=b_2=\f{2c^2(\rho_0)}{u_0^4},\ \ a_2=\f{(\gamma-1)c^2(\rho_0)}{u_0^2}+\f{2c^4(\rho_0)}{u_0^4},\\
b_1=\f{2}{u_0^4},
d_1=-a_2\p_1 s_0+b_2\p_1\fai_0=-\f{\gamma c^2(\rho_0)E_0}{(c^2(\rho_0)-u_0^2)u_0^2},\\
d_2=d_3=(\f{1}{u_0^2}\p_1\fai_0-\f{c^2(\rho_0)}{u_0^2}\p_1s_0)=-\f{E_0}{(c^2(\rho_0)-u_0^2)},\\
d_4=d_5=a_1\p_1s_0-b_1\p_1\fai_0=\f{2E_0}{(c^2(\rho_0)-u_0^2)u_0^2}.
\ea\right.\ee

We emphasize that $F_i(({\bf  W},\nabla {\bf  W}),i=1,2,3,5$ are second order terms of $W_i, i=1,\cdots,5$ and does not contain any space derivative of $W_2,W_3,W_4$.

Then it is easy to derive the elliptic system satisfied by $W_1$ and $W_5$:
\be\label{ee25} \left\{\begin{array}{ll}
\p_1\bigg((\f{c^2(\rho_0)}{u_0^2}-1)\p_1W_1-d_1W_1-d_5W_5\bigg)+\p_2\bigg((\f{c^2(\rho_0)}{u_0^2}-1)W_2\p_1W_1+\f{c^2(\rho_0)}{u_0^2}\p_2 W_1\bigg)\\
+\p_3\bigg((\f{c^2(\rho_0)}{u_0^2}-1)W_3\p_1W_1+\f{c^2(\rho_0)}{u_0^2}\p_3 W_1\bigg)+d_2\bigg((\f{c^2(\rho_0)}{u_0^2}-1)\p_1W_1-d_1W_1-d_5W_5\bigg)\\
-\f{1}{u_0^2}\rho_0 W_1-\bigg(\f{d_2}{u_0^2}+\p_1(\f{1}{u_0^2})\bigg)\p_1W_5=G_1({\bf  W},\nabla {\bf  W}),\\
\Delta W_5=\rho_0 W_1-\epsilon\t b(x)+F_5({\bf  W},\nabla {\bf  W}).\\
\bigg((\f{c^2(\rho_0)}{u_0^2}-1)\p_1W_1-d_1W_1-d_5W_5\bigg)(0,x_2,x_3)=F_1({\bf  W},\nabla {\bf  W})+\epsilon\bigg(\p_2\beta_2^{in}+\p_3\beta_3^{in}+\f{E^{in}}{u_0^2}\bigg),\\
\p_1 W_5(0,x_2,x_3)=\epsilon E^{in}(x_2,x_3),\\
W_1(1,x_2,x_3)=\epsilon s_e(x_2,x_3),\\
W_5(1,x_2,x_3)=0.
\end{array}\right. \ee
Here
\be \begin{aligned}
G_1({\bf  W},\nabla {\bf  W})=&\p_1 (F_1({\bf  W},\nabla {\bf  W})+d_4W_4)+\p_2 F_2({\bf  W},\nabla {\bf  W})+\p_3 F_3({\bf  W},\nabla {\bf  W})
\\&+\p_2\bigg(d_1W_1W_2+d_5W_2W_5+\f{1}{u_0^2}W_2\p_1W_5+W_2F_1({\bf  W},\nabla {\bf  W})+d_4W_2W_4\bigg)\\&+\p_3\bigg(d_1W_1W_3+d_5W_3W_5+\f{1}{u_0^2}W_3\p_1W_5+W_3F_1({\bf  W},\nabla {\bf  W})+d_5W_2W_5\bigg)
\\&-2(\p_2W_3\p_3W_2-\p_2W_2\p_3W_3)+d_2d_4W_4+d_2F_1({\bf  W},\nabla {\bf  W})\\&+\f{1}{u_0^2}\bigg[-\epsilon\t b(x)+F_5({\bf  W},\nabla {\bf  W})\bigg].
\end{aligned}
\ee
And $W_2, W_3$ and $W_4$ will satisfy the following hyperbolic equations respectively:
\be\label{ee27} \left\{\begin{array}{ll}
\p_1W_2+W_2\p_2W_2+W_3\p_3W_2+d_2W_2+\f{c^2(\rho_0)}{u_0^2}\p_2 W_1-\f{1}{u_0^2}\p_2 W_5=F_2({\bf  W},\nabla {\bf  W}),\\
\p_1W_3+W_2\p_2W_3+W_3\p_3W_3+d_3W_3+\f{c^2(\rho_0)}{u_0^2}\p_3 W_1-\f{1}{u_0^2}\p_3 W_5=F_3({\bf  W},\nabla {\bf  W}),\\
W_2(0,x_2,x_3)=\epsilon \beta_{20}(x_2,x_3),\\
W_3(0,x_2,x_3)=\epsilon \beta_{30}(x_2,x_3).
\end{array}\right. \ee

\be\label{ee28} \left\{\begin{array}{ll}
\p_1W_4+W_2\p_2W_4+W_3\p_3W_4=0,\\
 W_4(0,x_2,x_3)=\epsilon B^{in}(x_2,x_3).
\end{array}\right. \ee

\subsection{Proof of Theorem \ref{EPTH}}\label{EulerPoissonEPTH}\hspace*{\parindent}

The main idea is simple: we construct an operator $\Lambda$ on a suitable space, which will be bounded in a high order norm and contraction in a low order norm.

The solution class will be given by
$$\Xi=\bigg\{{\bf W}=(W_1,W_2,W_3,W_4,W_5): \|{\bf W}\|_{\Xi}=\|W_1,W_5\|_{C^{3,\alpha}(\O_e)}+\sum_{i=2}^4\|W_i\|_{C^{2,\alpha}(\O_e)}\leq \delta.\bigg\}$$
Here $\delta$ will be determined later. For a given ${\bf \t W}=(\t{W_1},\t{W_2},\t {W_3},\t {W_4},\t {W_5})\in \Xi$, we define an operator $\Lambda:{\bf \t W}=(\t{W_1},\t{W_2},\t {W_3},\t {W_4},\t {W_5})\mapsto {\bf W}=(W_1,W_2, W_3, W_4, W_5)$ mapping from $\Xi$ to itself, through the following iteration scheme.

Step 1. To resolve $W_4$.
\be\label{ee151} \left\{\begin{array}{ll}
\p_1W_4+\t W_2\p_2W_4+\t W_3\p_3W_4=0,\\
 W_4(0,x_2,x_3)=\epsilon B^{in}(x_2,x_3).
\end{array}\right. \ee
The particle path $(\tau,\t x_2(\tau;{\bf x}), \t x_3(\tau;{\bf x}))$ through $(x_1,x_2,x_3)$, is defined by the following ordinary differential equations:
\be\label{ee12} \left\{\begin{array}{ll}
\f{d \t x_2(\tau;{\bf x})}{d\tau}=\t W_2(\tau,\t x_2(\tau;{\bf x}),\t x_3(\tau;{\bf x})),\\
\f{d \t x_3(\tau;{\bf x})}{d\tau}=\t W_3(\tau,\t x_2(\tau;{\bf x}),\t x_3(\tau;{\bf x})),\\
\t x_2(x_1;{\bf x})=x_2,\\
\t x_2(x_1;{\bf x})=x_3.
\end{array}\right. \ee
Since $(\t W_2,\t W_3)\in C^{2,\alpha}(\O_e)$, it is easy to prove that $(\t x_2(0;{\bf x}),\t x_3(0;{\bf x}))$ belong to $C^{2,\alpha}(\O_e)$ with respect to ${\bf x}$. Indeed, the following estimate holds:
\be\label{ee13}
\|(\t x_2(0;{\bf x}),\t x_3(0;{\bf x}))\|_{C^{2,\alpha}(\O_e)}\leq C_1.
\ee

One can find a unique solution $ W_4(x)=\epsilon B^{in}(\t x_2(0;{\bf x}),\t x_3(0;{\bf x}))\in C^{1,\alpha}(\O_e)$ to (\ref{ee151}). Furthermore, the following estimate holds:
\be\label{ee152}
\|W_4\|_{C^{2,\alpha}(\O_e)}\leq \epsilon \|B^{in}\|_{C^{2,\alpha}(\mathrm{T}^2)}\leq C_2\epsilon.
\ee

Step 2. To obtain $W_1$ and $W_5$ by solving the following linearized second order elliptic systems:
\be\label{ee9} \left\{\begin{array}{ll}
\p_1\bigg((\f{c^2(\rho_0)}{u_0^2}-1)\p_1W_1-d_1W_1-d_5W_5\bigg)+\p_2\bigg((\f{c^2(\rho_0)}{u_0^2}-1)\t W_2\p_1W_1+\f{c^2(\rho_0)}{u_0^2}\p_2 W_1\bigg)\\
+\p_3\bigg((\f{c^2(\rho_0)}{u_0^2}-1)\t W_3\p_1W_1+\f{c^2(\rho_0)}{u_0^2}\p_3 W_1\bigg)+d_2\bigg((\f{c^2(\rho_0)}{u_0^2}-1)\p_1W_1-d_1W_1-d_5W_5\bigg)\\
-\f{1}{u_0^2}\rho_0 W_1-\bigg(\f{d_2}{u_0^2}+\p_1(\f{1}{u_0^2})\bigg)\p_1W_5=G_1({\bf \t  W},\nabla {\bf\t  W}),\\
\Delta W_5=\rho_0 W_1-\epsilon\t b(x)+F_5({\bf \t W},\nabla {\bf \t  W}).\\
\bigg((\f{c^2(\rho_0)}{u_0^2}-1)\p_1W_1-d_1W_1-d_5W_5\bigg)(0,x_2,x_3)=F_1({\bf\t  W},\nabla {\bf \t W})+d_4W_4\\ \quad \quad\quad\quad+\epsilon\bigg(\p_2\beta_2^{in}+\p_3\beta_3^{in}+\f{1}{u_0^2}E^{in}\bigg),\\
\p_1 W_5(0,x_2,x_3)=\epsilon E^{in}(x_2,x_3),\\
W_1(1,x_2,x_3)=\epsilon s_e(x_2,x_3),\\
W_5(1,x_2,x_3)=0.
\end{array}\right. \ee

Take $U_1=W_1-\epsilon s_e(x_2,x_3), U_5=W_5$, then $U_1,U_5$ satisfy the following elliptic systems:
\be\label{ee91} \left\{\begin{array}{ll}
\p_1\bigg((\f{c^2(\rho_0)}{u_0^2}-1)\p_1U_1-d_1U_1-d_5U_5\bigg)+\p_2\bigg((\f{c^2(\rho_0)}{u_0^2}-1)\t W_2\p_1U_1+\f{c^2(\rho_0)}{u_0^2}\p_2 U_1\bigg)\\
+\p_3\bigg((\f{c^2(\rho_0)}{u_0^2}-1)\t W_3\p_1U_1+\f{c^2(\rho_0)}{u_0^2}\p_3 U_1\bigg)+d_2\bigg((\f{c^2(\rho_0)}{u_0^2}-1)\p_1U_1-d_1U_1-d_5U_5\bigg)\\
-\f{1}{u_0^2}\rho_0 U_1-\bigg(\f{d_2}{u_0^2}+\p_1(\f{1}{u_0^2})\bigg)\p_1U_5=\t G_1({\bf \t  W},\nabla {\bf\t  W}),\\
\Delta U_5=\rho_0 U_1+\t G_5({\bf \t  W},\nabla {\bf\t  W}).\\
\bigg((\f{c^2(\rho_0)}{u_0^2}-1)\p_1U_1-d_1U_1-d_5U_5\bigg)(0,x_2,x_3)=\t H_1({\bf \t  W},\nabla {\bf\t  W}),\\
\p_1 U_5(0,x_2,x_3)=\t H_2({\bf \t  W},\nabla {\bf\t  W})=\epsilon E^{in}(x_2,x_3),\\
U_1(1,x_2,x_3)=0,\\
U_5(1,x_2,x_3)=0.
\end{array}\right. \ee
Here
\be\begin{cases}\begin{aligned}
\t G_1({\bf \t  W},\nabla {\bf\t  W})&=G_1({\bf \t  W},\nabla {\bf\t  W})+\epsilon\p_1(d_1s_e)
-\sum_{i=2}^3\epsilon\p_i(\f{c^2(\rho_0)}{u_0^2}\p_is_e)+\epsilon (d_1d_2s_e+\f{\rho_0}{u_0^2}s_e).\\
\t G_5({\bf \t  W},\nabla {\bf\t  W})&=-\epsilon\t b(x)+F_5({\bf \t W},\nabla {\bf \t  W})+\epsilon \rho_0s_e,\\
\t H_1({\bf \t  W},\nabla {\bf\t  W})&=F_1({\bf\t  W},\nabla {\bf \t W})+d_4W_4+\epsilon\bigg(\p_2\beta_2^{in}+\p_3\beta_3^{in}+\f{1}{u_0^2}E^{in}\bigg)+\epsilon d_1s_e.
\end{aligned}\end{cases}\ee
Set $\mathbb{H}=\{U(x)\in H^1(\O_e):U(1,x_2,x_3)=0.\}$. The weak formulation for (\ref{ee91}) is the following: To find ${\bf U}=(U_1,U_5)\in \mathbb{H}\times\mathbb{H}$,
such that for any ${\bf V}=(V_1,V_5)\in \mathbb{H}\times\mathbb{H}$, we have
$\mathbf{B}({\bf U},{\bf V})=\mathfrak{F}({\bf V})$.
Here
\be\begin{aligned}\label{ee92}
\mathbf{B}({\bf U},{\bf V})=&\int_{\O_e}\bigg\{\bigg((\f{c^2(\rho_0)}{u_0^2}-1)\p_1U_1-d_1U_1-d_5U_5\bigg)\p_1V_1\\&+\bigg((\f{c^2(\rho_0)}{u_0^2}-1)\t W_2\p_1U_1+\f{c^2(\rho_0)}{u_0^2}\p_2 U_1\bigg)\p_2V_1\\&
+\bigg((\f{c^2(\rho_0)}{u_0^2}-1)\t W_3\p_1U_1+\f{c^2(\rho_0)}{u_0^2}\p_3 U_1\bigg)\p_3V_1\bigg\}dx\\&
-\int_{\O_e}d_2\bigg((\f{c^2(\rho_0)}{u_0^2}-1)\p_1U_1-d_1U_1-d_5U_5\bigg)V_1dx\\&
+\int_{\O_e}\bigg\{\f{\rho_0}{u_0^2} U_1V_1+\bigg(\f{d_2}{u_0^2}+\p_1(\f{1}{u_0^2})\bigg)\p_1U_5V_1\bigg\}dx\\&+\int_{\O_e}\nabla U_5\cdot\nabla V_5dx+\int_{\O_e}\rho_0U_1V_5dx.
\end{aligned}\ee
\be\begin{aligned}\label{ee93}
\mathfrak{F}({\bf V})=\int_{\O_e}\t G_1 V_1dx-\int_{\O_e}\t G_2 V_5dx-\int_{\mathrm{T}^2}(\t H^1V_1+\t H^2V_5)(0,x_2,x_3)dx_2dx_3.
\end{aligned}\ee
A direct computation shows that $\mathbf{B}$ is a bounded bilinear operator in $\mathbb{H}\times\mathbb{H}$ and $\mathfrak{F}$ is a bounded linear operator in $\mathbb{H}\times\mathbb{H}$. Also $\mathbf{B}$ is coercive. Hence same arguments as in Chapter 8 in \cite{GT} shows us that
Fredholm alternative theorem holds for (\ref{ee91}). Hence uniqueness theorem will imply the existence of solution to (\ref{ee91}).
 It suffices to show that $\mathbf{B}({\bf U},{\bf U})=0$ implies ${\bf U}=0$.
\be\begin{aligned}\label{ee94}
\mathbf{B}({\bf U},{\bf U})=&\int_{\O_e}\bigg\{(\f{c^2(\rho_0)}{u_0^2}-1)(\p_1U_1)^2+\f{c^2(\rho_0)}{u_0^2}[(\p_2U_1)^2+(\p_3U_1)^2+|\nabla U_5|^2]\bigg\}dx\\&
+\int_{\O_e}\bigg\{(\f{c^2(\rho_0)}{u_0^2}-1)(\t W_2\p_1U_1\p_2U_1+\t W_3\p_1U_1\p_3U_1)\bigg\}dx\\&-\int_{\O_e}[d_1+d_2(\f{c^2(\rho_0)}{u_0^2}-1)]U_1\p_1U_1dx-\int_{\O_e}d_5\p_1U_1U_5dx\\&
+\int_{\O_e}\bigg(\f{d_2}{u_0^2}+\p_1(\f{1}{u_0^2})\bigg)\p_1U_5U_1dx+\int_{\O_e}\bigg(\f{\rho_0}{u_0^2}+d_1d_2\bigg)U_1^2dx\\&
+\int_{\O_e}\bigg(\rho_0+d_2d_5\bigg)U_1U_5dx=I+\sum_{i=1}^6 I_i=0.
\end{aligned}\ee
Now we estimate these terms respectively.
\be\begin{aligned}
I_6&\leq \f{1}{4\eta}\int_{\O_e}(\rho_0+d_2d_5)^2U_1^2dx+\eta\|U_5\|_{L^2(\O_e)}^2\\&\leq \f{1}{4\eta}\int_{\O_e}(\rho_0+d_2d_5)^2U_1^2dx+\eta C_{\O_e}\|\nabla U_5\|_{L^2(\O_e)}^2\\&
\leq \f{C_{\O_e}}{2}\bigg(\max_{x\in[0,1]}\rho_0^2+C(b_0)\max_{x\in[0,1]}(\f{E_0}{u_0^2})\bigg)\|U_1\|_{L^2(\O_e)}^2+\f{1}{2}\|\nabla U_5\|_{L^2(\O_e)}^2,
\end{aligned}\ee
\be\begin{aligned}
I_5&\geq\min_{x\in[0,1]}(\f{\rho_0}{u_0^2})\|U_1\|_{L^2(\O_e)}^2-\max_{x\in[0,1]}(\f{E_0^2}{u_0^2})\|U_1\|_{L^2(\O_e)}^2\\&
\geq \f{\min_{x\in[0,1]}\rho_0}{\max_{x\in[0,1]}u_0^2}\|U_1\|_{L^2(\O_e)}^2-\max_{x\in[0,1]}(\f{E_0^2}{u_0^2})\|U_1\|_{L^2(\O_e)}^2.
\end{aligned}\ee
\be\begin{aligned}
\sum_{i=1}^4|I_i|&\leq C(b_0)\bigg(\max_{x\in[0,1]}(\f{1}{u_0^2})\delta\|\nabla U_1\|_{L^2(\O_e)}^2\\&+\max_{x\in[0,1]}(\f{E_0}{u_0^2})[\| U_1\|_{L^2(\O_e)}^2
+\| U_5\|_{L^2(\O_e)}^2+\|\nabla U_1\|_{L^2(\O_e)}^2+\|\nabla U_5\|_{L^2(\O_e)}^2]\bigg).
\end{aligned}\ee
Hence we have
\be\begin{aligned}
0=I+\sum_{i=1}^6 I_i&\geq \int_{\O_e}|\nabla U_1|^2+|\nabla U_5|^2dx-C(b)\bigg(\max_{x\in[0,1]}(\f{1}{u_0^2})\delta\|\nabla U_1\|_{L^2(\O_e)}^2\\&+\max_{x\in[0,1]}(\f{E_0}{u_0^2})[\| U_1\|_{L^2(\O_e)}^2
+\| U_5\|_{L^2(\O_e)}^2+\|\nabla U_1\|_{L^2(\O_e)}^2\\&
+\|\nabla U_5\|_{L^2(\O_e)}^2]\bigg)+\bigg(\f{\min_{x\in[0,1]}\rho_0}{\max_{x\in[0,1]}u_0^2}-\f{C_{\O_e}}{2}\max_{x\in[0,1]}\rho_0^2\bigg)\|U_1\|_{L^2(\O_e)}^2
\\&-\max_{x\in[0,1]}(\f{E_0^2}{u_0^2})\|U_1\|_{L^2(\O_e)}^2
-\f{C_{\O_e}}{2}C(b_0)\max_{x\in[0,1]}(\f{E_0}{u_0^2})\|U_1\|_{L^2(\O_e)}^2-\f{1}{2}\|\nabla U_5\|_{L^2(\O_e)}^2\\&
\geq \f{1}{2} \int_{\O_e}|\nabla U_1|^2+|\nabla U_5|^2dx-C(b_0)C_{\O_e}\bigg(\max_{x\in[0,1]}(\f{1}{u_0^2})\delta\|\nabla U_1\|_{L^2(\O_e)}^2\\&+\max_{x\in[0,1]}(\f{E_0}{u_0^2})[\| U_1\|_{L^2(\O_e)}^2
+\| U_5\|_{L^2(\O_e)}^2+\|\nabla U_1\|_{L^2(\O_e)}^2
+\|\nabla U_5\|_{L^2(\O_e)}^2]\bigg)\\&\geq \f{1}{4}(\|\nabla U_1\|_{L^2(\O_e)}^2
+\|\nabla U_5\|_{L^2(\O_e)}^2).
\end{aligned}\ee
Here we have used some special properties (\ref{eb25}) of our background solutions.

This implies that $U_1=U_5=0$. We establish the uniqueness for (\ref{ee91}). By Fredholm alternative theorem, there exists a unique solution $(W_1,W_5)$ to (\ref{ee9}).

By standard elliptic estimates in \cite{ADN}, (\ref{ee9}) has a unique solution $(W_1,W_5)\in C^{3,\alpha}(\bar\O)\times C^{3,\alpha}(\bar\O)$ and satisfies the following estimate:
\be\begin{aligned}\label{ee10}
\|W_1,W_5\|_{C^{3,\alpha}(\O_e)}&\leq C_3\bigg(\|G_1({\bf \t W},\nabla {\bf \t W})\|_{C^{1,\alpha}(\O_e)}+\epsilon\|\t b(x)\|_{C^{1,\alpha}(\O_e)}+\|F_5({\bf \t W},\nabla {\bf \t W})\|_{C^{1,\alpha}(\O_e)}\\&
+\epsilon\|\p_2\beta_{2}^{in}+\p_3\beta_{3}^{in}+\f{1}{u_0^2}E^{in}\|_{C^{2,\alpha}(\mathrm{T}^2)}+\|F_1({\bf \t W},\nabla {\bf \t W})+d_4W_4\|_{C^{2,\alpha}(\mathrm{T}^2)}
\\&+\epsilon\|E^{in}\|_{C^{2,\alpha}(\mathrm{T}^2)}+\epsilon\|s_e\|_{C^{3,\alpha}(\mathrm{T}^2)}\bigg)\\&\leq C_4(\delta^2+\epsilon).
\end{aligned}\ee
Here $C_3$ depends only on background solution and $\O_e$, $C_4$ depends also on $\beta_{2}^{in},\beta_{3}^{in},B^{in}$ and $E^{in},\t b, s_e$.

Step 3. To obtain $W_2,W_3$ by solving the following hyperbolic equations:
\be\label{ee11} \left\{\begin{array}{ll}
\p_1W_2+\t W_2\p_2W_2+\t W_3\p_3W_2+d_2W_2+\f{c^2(\rho_0)}{u_0^2}\p_2 W_1-\f{1}{u_0^2}\p_2 W_5=F_2({\bf \t W},\nabla {\bf \t W}),\\
\p_1W_3+\t W_2\p_2W_3+\t W_3\p_3W_3+d_3W_3+\f{c^2(\rho_0)}{u_0^2}\p_3 W_1-\f{1}{u_0^2}\p_3 W_5=F_3({\bf \t W},\nabla {\bf \t W}),\\
W_2(0,x_2,x_3)=\epsilon \beta_{2}^{in}(x_2,x_3),\\
W_3(0,x_2,x_3)=\epsilon \beta_{3}^{in}(x_2,x_3).
\end{array}\right. \ee

Then by the characteristic methods, we have the following formulas:
\be\label{ee14} \left\{\begin{array}{ll}
W_2({\bf x})=&\epsilon e^{-\int_0^{x_1}d_2(s)ds}\beta_{2}^{in}(\t x_2(0;{\bf x}),\t x_3(0;{\bf x}))+\int_0^{x_1}e^{-\int_{\tau}^{x_1}d_2(s)ds}\bigg(F_2({\bf \t W},\nabla {\bf \t W})\\&-\f{c^2(\rho_0)}{u_0^2}\p_2 W_1+\f{1}{u_0^2}\p_2 W_5\bigg)(\tau,\t x_2(\tau;{\bf x}),\t x_3(\tau;{\bf x}))d\tau,\\
W_3({\bf x})=&\epsilon e^{-\int_0^{x_1}d_3(s)ds}\beta_{3}^{in}(\t x_2(0;{\bf x}),\t x_3(0;{\bf x}))+\int_0^{x_1}e^{-\int_{\tau}^{x_1}d_3(s)ds}\bigg(F_3({\bf \t W},\nabla {\bf \t W})\\&-\f{c^2(\rho_0)}{u_0^2}\p_3 W_1+\f{1}{u_0^2}\p_3 W_5\bigg)(\tau,\t x_2(\tau;{\bf x}),\t x_3(\tau;{\bf x}))d\tau.
\end{array}\right. \ee
These enable one to obtain the following estimate:
\be\begin{aligned}\label{ee15}
\|(W_2({\bf x}),W_3({\bf x}))\|_{C^{2,\alpha}(\O_e)}&\leq C_5\bigg[\epsilon \|(\beta_{2}^{in},\beta_{3}^{in})\|_{C^{2,\alpha}(\mathrm{T}^2)}
+\|{\bf \t W}\|_{\Xi}^2+\|W_1,W_5\|_{C^{3,\alpha}(\O_e)}\bigg]\\&\leq C_6(\delta^2+\epsilon).
\end{aligned}\ee

This, together with the estimates (\ref{ee10}), (\ref{ee15}) and (\ref{ee152}), gives
\be\label{ee153}
\|{\bf W}\|_{\Xi}=\|W_1,W_5\|_{C^{3,\alpha}(\O_e)}+\|W_2,W_3,W_4\|_{C^{2,\alpha}(\O_e)}\leq C_7(\delta^2+\epsilon).
\ee
Here $C_7=max\{C_2,C_4,C_6\}$, which depends only on background solution and $\beta_{2}^{in},\beta_{3}^{in},B^{in}$ and $E^{in},\t b, s_e$.

Choose $\epsilon_1$ small enough, such that if $\epsilon\leq \epsilon_1$, then $C_7^2\epsilon <\f{1}{4}$. Set $\delta=2C_7\epsilon$, then $C_7\delta<\f{1}{2}$ and $C_7(\delta^2+\epsilon)<C_7\epsilon+\f{1}{2}\delta=\delta$. This implies that $\Lambda$ maps $\Xi$ to itself.

It remains to show that the mapping $\Lambda: \Xi\rightarrow \Xi$ is a contraction operator.
Suppose $\Lambda:(\t W_1^k,\t W_2^k, \t W_3^k,\t W_4^k,\t W_5^k)\mapsto (W_1^k, W_2^k, W_3^k, W_4^k, W_5^k)$ for $k=1,2$.
 Define the difference $(Y_1,Y_2, Y_3, Y_4)=(W_1^1-W_1^2,W_2^1-W_2^2,W_3^1-W_3^2, W_4^1-W_4^2,W_5^1-W_5^2)$
 and $(\t Y_1,\t Y_2,\t Y_3,\t Y_4)=(\t W_1^1-\t W_1^2,\t W_2^1-\t W_2^2,\t W_3^1-\t W_3^2,\t W_4^1-\t W_4^2,\t W_5^1-\t W_5^2)$.

Step 1. Estimate of $Y_4$.

Indeed, $Y_4$ satisfies the following equation:
\be\label{ee21} \left\{\begin{array}{ll}
\p_1\bar Y_4+\t W_2^1\p_2Y_4+\t W_3^1\p_3Y_4=-\t Y_2\p_2 W_4^2-\t Y_3\p_3 W_4^2,\\
 \bar Y_4(0,x_2,x_3)=0.
\end{array}\right. \ee
Then the following estimate holds:
\be\begin{aligned}\label{ee22}
\|Y_4\|_{C^{1,\alpha}(\O_e)}\leq C_{8}\delta \|(\t Y_2, \t Y_3)\|_{C^{1,\alpha}(\O_e)}.
\end{aligned}\ee

Step 2. Estimates of $Y_1$ and $Y_5$.

$Y_1$ satisfies the following elliptic system:

\be\label{ee16} \left\{\begin{array}{ll}
\p_1\bigg((\f{c^2(\rho_0)}{u_0^2}-1)\p_1Y_1-d_1Y_1-d_5Y_5\bigg)+\p_2\bigg((\f{c^2(\rho_0)}{u_0^2}-1)\t W_2^1\p_1Y_1+\f{c^2(\rho_0)}{u_0^2}\p_2 Y_1\bigg)\\
+\p_3\bigg((\f{c^2(\rho_0)}{u_0^2}-1)\t W_3^1\p_1Y_1+\f{c^2(\rho_0)}{u_0^2}\p_3 Y_1\bigg)+d_2\bigg((\f{c^2(\rho_0)}{u_0^2}-1)\p_1Y_1-d_1Y_1-d_5Y_5\bigg)\\
-\f{1}{u_0^2}\rho_0 Y_1-\bigg(\f{d_2}{u_0^2}+\p_1(\f{1}{u_0^2})\bigg)\p_1Y_5=H_1,\\
\Delta Y_5=\rho_0 Y_1+F_5({\bf \t W^1},\nabla {\bf \t W^1})-F_5({\bf \t W^2},\nabla {\bf \t W^2}),\\
\bigg((\f{c^2(\rho_0)}{u_0^2}-1)\p_1Y_1-d_1Y_1-d_5Y_5\bigg)(0,x_2,x_3)=F_1({\bf \t W^1},\nabla {\bf \t W^1})-F_1({\bf \t W^2},\nabla {\bf \t W^2}),\\
\p_1Y_5(0,x_2,x_3)=0,\\
Y_1(1,x_2,x_3)=0,\\
Y_5(1,x_2,x_3)=0.
\end{array}\right. \ee

$$H_1=G_1({\bf \t W^1},\nabla {\bf \t W^1})-G_1({\bf \t W^2},\nabla {\bf \t W^2})-\p_2\bigg((\f{c^2(\rho_0)}{u_0^2}-1)\t Y_2\p_1W_1^2\bigg)-\p_3\bigg((\f{c^2(\rho_0)}{u_0^2}-1)\t Y_3\p_1W_1^2\bigg).$$

Then the following estimate holds:
\be\begin{aligned}\label{ee18}
\|(Y_1,Y_5)\|_{C^{2,\alpha}(\O_e)}&\leq C_9\bigg(\|H_1\|_{C^{\alpha}(\O_e)}+\|F_5({\bf \t W^1},\nabla {\bf \t W^1})-F_5({\bf \t W^2},\nabla {\bf \t W^2})\|_{C^{1,\alpha}(\O_e)}\\&
+\|F_1({\bf \t W^1},\nabla {\bf \t W^1})-F_1({\bf \t W^2},\nabla {\bf \t W^2})+d_4Y_4\|_{C^{1,\alpha}(\mathrm{T}^2)}\bigg)
\\ &\leq C_{10}\delta \bigg(\|\t Y_1,\t Y_5\|_{C^{2,\alpha}(\O_e)}+\|\t Y_2,\t Y_3,\t Y_4\|_{C^{1,\alpha}(\O_e)}\bigg).
\end{aligned}\ee

Step 3. Estimate of $Y_2,Y_3$.

It follows from (\ref{ee11}), $Y_2$ and $Y_3$ satisfy the following system:
\be\label{ee19} \left\{\begin{array}{ll}
\p_1 Y_2+\t W_2^1\p_2Y_2+\t W_3^1\p_3Y_2+d_2Y_2+\f{c^2(\rho_0)}{u_0^2}\p_2Y_1-\f{1}{u_0^2}\p_2Y_5=K_1,\\
\p_1 Y_3+\t W_2^1\p_2Y_3+\t W_3^1\p_3Y_3+d_3Y_3+\f{c^2(\rho_0)}{u_0^2}\p_3Y_1-\f{1}{u_0^2}\p_3Y_5=K_2.\\
Y_2(0,x_2,x_3)=0,\\
Y_3(0,x_2,x_3)=0.
\end{array}\right. \ee
Here

\be \left\{\begin{array}{ll}
K_1=-(\t Y_2\p_2W_2^2+\t Y_3\p_3W_2^2)+F_2({\bf \t W^1},\nabla {\bf \t W^1})-F_2({\bf \t W^2},\nabla {\bf \t W^2}),\\
K_2=-(\t Y_2\p_2W_3^2+\t Y_3\p_3W_3^2)+F_3({\bf \t W^1},\nabla {\bf \t W^1})-F_3({\bf \t W^2},\nabla {\bf \t W^2}).
\end{array}\right.\ee

Then the following estimate holds:
\be\begin{aligned}\label{ee20}
\|(Y_2, Y_3)\|_{C^{1,\alpha}(\O_e)}&\leq C_{11}\bigg(\delta \|(\t Y_2, \t Y_3, \t Y_4)\|_{C^{1,\alpha}(\O_e)}+\delta\|\nabla Y_1,\nabla Y_5\|_{C^{1,\alpha}(\O_e)}\bigg)
\\&\leq C_{12}\delta\bigg(\|\t Y_1,\t Y_5\|_{C^{2,\alpha}(\O_e)}+\|\t Y_2,\t Y_3,\t Y_4\|_{C^{1,\alpha}(\O_e)}\bigg).
\end{aligned}\ee

Setting $C_{13}=max\{C_8,C_{10},C_{12}\}$, then take $\epsilon_2$ small enough such that $C_{13}\epsilon_2<1$. Now we
 take $\epsilon_0=min\{\epsilon_1,\epsilon_2\}$, then if $0<\epsilon<\epsilon_0$, $\Lambda$ maps $\Xi$ to itself a
 nd is contraction in low order norm. Hence $\Lambda$ has a unique fixed point, which will be the solution to (\ref{ea1}). We have finished our proof.

\begin{remark}{\it
We can not generalize this result to a general 3-D nozzle. In general, the velocity field $\beta_2$ and $\beta_3$ will lose one order derivative when integrating along the particle path. The delicate nonlinear coupling between the hyperbolic modes and elliptic modes $(\beta_2,\beta_3)$ makes it extremely difficult to develop an effective iteration scheme in a general 3-D nozzle. How to effectively explore the good property of the quantity $W=\p_2\beta_3-\p_3\beta_2+\beta_3\p_1\beta_2-\beta_2\p_1\beta_3$ will be investigated in the forthcoming paper.
}
\end{remark}

\section {Appendix}\label{appendix}\hspace*{\parindent}

In this appendix, we give the detailed calculations for (\ref{eb26}). First we rewrite (\ref{eb25}) as follows:
\be\label{appendix1}\left\{\ba{l}
G^{-1}{\bf D}\beta_2+\f{1}{2(B+\varphi-h(\rho))}(\beta_2\p_1-\p_2)(\varphi-h(\rho))=0,\\
G^{-1}{\bf D}\beta_3+\f{1}{2(B+\varphi-h(\rho))}(\beta_3\p_1-\p_3)(\varphi-h(\rho))=0.
\ea\right.
\ee

Applying $\beta_3\p_1-\p_3$ and $-(\beta_2\p_1-\p_2)$ to (\ref{appendix1}) and adding them together, we obtain

\be\begin{aligned}
0&=(\beta_3\p_1-\p_3)(G^{-1}{\bf D}\beta_2)-(\beta_2\p_1-\p_2)(G^{-1}{\bf D}\beta_3)\\&\quad\quad+\f{1}{2(B+\varphi-h(\rho))}\bigg[(\beta_3\p_1-\p_3)(\beta_2\p_1-\p_2)-(\beta_2\p_1-\p_2)(\beta_3\p_1-\p_3)\bigg](\varphi-h(\rho))\\&
\quad\quad-\f{1}{2(B+\varphi-h(\rho))^2} \bigg[(\beta_3\p_1-\p_3) B (\beta_2\p_1-\p_2)(\varphi-h(\rho))\\&
\quad\quad\quad\quad-(\beta_2\p_1-\p_2) B (\beta_3\p_1-\p_3)(\varphi-h(\rho)) \bigg]\\&
 =(\beta_3\p_1-\p_3)(G^{-1}{\bf D}\beta_2)-(\beta_2\p_1-\p_2)(G^{-1}{\bf D}\beta_3)+\f{1}{2(B+\varphi-h(\rho))} \p_1 (\varphi-h(\rho))W
\\& \quad  -\f{1}{2(B+\varphi-h(\rho))^2} \bigg[(1,\beta_2,\beta_3)\cdot \bigg(\nabla (\varphi-h(\rho))\times \nabla B \bigg)\bigg].
\end{aligned}\ee

While
\be
\begin{aligned}
J:&=(\beta_3\p_1-\p_3)(G^{-1}{\bf D}\beta_2)-(\beta_2\p_1-\p_2)(G^{-1}{\bf D}\beta_3)\\
&=G^{-1}{\bf D}((\beta_3\p_1-\p_3)\beta_2-(\beta_2\p_1-\p_2)\beta_3)\\&+G^{-1}\sum_{j=1}^3((\beta_3\p_1-\p_3)\beta_j\p_j\beta_2-(\beta_2\p_1-\p_2)\beta_j\p_j\beta_3)\\&
   +G^{-1}\sum_{j=1}^3\beta_j\bigg\{[(\beta_3\p_1-\p_3)\p_j-\p_j(\beta_3\p_1-\p_3)]\beta_2
   \\&\quad\quad\quad\quad-[(\beta_2\p_1-\p_2)\p_j-\p_j(\beta_2\p_1-\p_2)]\beta_3\bigg\}\\&
   +(\beta_3\p_1-\p_3)G^{-1}{\bf D}\beta_2-(\beta_2\p_1-\p_2)G^{-1}{\bf D}\beta_3\\ = &G^{-1}{\bf D}W+J_1+J_2+J_3.
\end{aligned}\ee

Now we compute $J_i, i=1,2,3$ respectively.

\be\begin{aligned}
J_1&=G^{-1}\bigg[W(\p_2\beta_2+\p_3\beta_3)-(\beta_3\p_1-\p_3)\beta_2\p_3\beta_3+(\beta_3\p_1-\p_3)\beta_3\p_3\beta_2
\\&\quad\quad\quad-(\beta_2\p_1-\p_2)\beta_2\p_2\beta_3+(\beta_2\p_1-\p_2)\beta_3\p_2\beta_2\bigg]\\&
= G^{-1}[W(\p_2\beta_2+\p_3\beta_3)+Z].
\end{aligned}\ee

Here $Z=\beta_3\p_1\beta_3\p_3\beta_2-\beta_3\p_1\beta_2\p_3\beta_3+\beta_2\p_1\beta_3\p_2\beta_2-\beta_2\p_1\beta_2\p_2\beta_3$.

\be\begin{aligned}
J_2&=G^{-1}\sum_{j=1}^3\beta_j\bigg\{[(\beta_3\p_1-\p_3)\p_j-\p_j(\beta_3\p_1-\p_3)]\beta_2\\&\quad\quad\quad\quad-[(\beta_2\p_1-\p_2)\p_j-\p_j(\beta_2\p_1-\p_2)]\beta_3\bigg\}\\&
    =G^{-1}Z.
\end{aligned}\ee

\be\begin{aligned}
J_3&=-G^{-2}\bigg[(\beta_3\p_1-\p_3)G{\bf D}\beta_2-(\beta_2\p_1-\p_2)G{\bf D}\beta_3\bigg]\\&=-2G^{-2}\sum_{j=1}^3\bigg[\beta_j(\beta_3\p_1-\p_3)\beta_j{\bf D}\beta_2-\beta_j(\beta_2\p_1-\p_2)\beta_j{\bf D}\beta_3\bigg]\\&
    = -2G^{-2}[W(\beta_2{\bf D}\beta_2+\beta_3{\bf D}\beta_3)+J_{31}].
\end{aligned}\ee
Here
\be\begin{aligned}
J_{31}&=\beta_2(\beta_2\p_1-\p_2)\beta_3{\bf D}\beta_2+\beta_3(\beta_3\p_1-\p_3)\beta_3{\bf D}\beta_2\\&\quad\quad-\beta_3(\beta_3\p_1-\p_3)\beta_2{\bf D}\beta_3-\beta_2(\beta_2\p_1-\p_2)\beta_2{\bf D}\beta_3\\&
       =\beta_2\bigg[(\beta_2\p_1-\p_2)\beta_3{\bf D}\beta_2-(\beta_2\p_1-\p_2)\beta_2{\bf D}\beta_3\bigg]\\&\quad\quad+\beta_3\bigg[(\beta_3\p_1-\p_3)\beta_3{\bf D}\beta_2-(\beta_2\p_1-\p_2)\beta_2{\bf D}\beta_3\bigg]\\&
       =\beta_2\bigg[(\beta_2\p_1-\p_2)\beta_3\p_1\beta_2-(\beta_2\p_1-\p_2)\beta_2\p_1\beta_3\bigg]
       \\&\quad+\beta_3\bigg[(\beta_3\p_1-\p_3)\beta_3\p_1\beta_2-(\beta_3\p_1-\p_3)\beta_2\p_1\beta_3\bigg]\\&
        \quad+\bigg[(\beta_2\p_1-\p_2)\beta_3(\beta_2^2\p_2\beta_2+\beta_2\beta_3\p_3\beta_2)-(\beta_2\p_1-\p_2)\beta_2(\beta_2^2\p_2\beta_3+\beta_2\beta_3\p_3\beta_3)\bigg]\\&
        \quad+\bigg[(\beta_3\p_1-\p_3)\beta_3(\beta_2\beta_3\p_2\beta_2+\beta_3^2\p_3\beta_2)-(\beta_3\p_1-\p_3)\beta_2(\beta_2\beta_3\p_2\beta_3+\beta_3^2\p_3\beta_3)\bigg]\\&  =GZ.
\end{aligned}\ee

Hence we have $J_3=-2G^{-2}[W(\beta_2 {\bf D}\beta_2+\beta_3 {\bf D}\beta_3)]-2G^{-1}Z$. Substitute these calculations into the above formula to get:

\be\begin{aligned}
0&=J+\f{1}{2(B+\varphi-h(\rho))} \p_1 (\varphi-h(\rho))W-\f{1}{2(B+\varphi-h(\rho))^2} \bigg[(1,\beta_2,\beta_3)\cdot \bigg(\nabla (\varphi-h(\rho))\times \nabla B \bigg)\bigg]\\&
 =G^{-1}{\bf D}W+G^{-1}\bigg(\p_2\beta_2+\p_3\beta_3+\f{1}{2(B+\varphi-h(\rho))}G\p_1(\varphi-h(\rho))\bigg)W-G^{-2}W {\bf D} G\\
 &\quad-\f{1}{2(B+\varphi-h(\rho))^2} \bigg[(1,\beta_2,\beta_3)\cdot \bigg(\nabla (\varphi-h(\rho))\times \nabla B \bigg)\bigg]\\&
 =G^{-1} {\bf D}W-G^{-1}W {\bf D}s-G^{-2}W {\bf D} G-\f{1}{2(B+\varphi-h(\rho))^2} \bigg[(1,\beta_2,\beta_3)\cdot \bigg(\nabla (\varphi-h(\rho))\times \nabla B \bigg)\bigg]\\
&=\rho {\bf D}(\f{W}{\rho G})-\f{1}{2(B+\varphi-h(\rho))^2} \bigg[(1,\beta_2,\beta_3)\cdot \bigg(\nabla (\varphi-h(\rho))\times \nabla B \bigg)\bigg].
\end{aligned}\ee

This implies that
\be
{\bf D}(\f{W}{\rho G})-\f{1}{2\rho(B+\varphi-h(\rho))^2} \bigg[(1,\beta_2,\beta_3)\cdot \bigg(\nabla (\varphi-h(\rho))\times \nabla B \bigg)\bigg]=0.
\ee

\bibliographystyle{plain}

\end{document}